\documentclass[12pt,twoside]{article}
\def\beq{\begin{equation}}
\def\eeq{\end{equation}}
\def\bea{\begin{eqnarray}}
\def\eea{\end{eqnarray}}
\def\ba{\begin{array}}
\def\ea{\end{array}}
\def\nn{\nonumber} 
\def\lrb{\left(}
\def\rrb{\right)} 
\def\lsb{\left[}
\def\rsb{\right]} 
\def\lcb{\left\{}
\def\rcb{\right\}}
\def\fr{\frac}
\def\ra{\rightarrow}
\def\lra{\longrightarrow}
\def\h{\hat} 
\def\lb{\label} 
\def\hsx{\hat{\cal{X}}}
\def\hsz{\hat{\cal{Z}}}
\def\hsa{\hat{\cal{A}}}
\def\hsn{\hat{\cal{N}}} 
\def\he{\hat{\eta}}
\begin{document}
\begin{center}
{\Large\bf  $(p,q)$-Rogers-Szeg\"{o} polynomial \\ 
and the $(p,q)$-oscillator}\footnote{{\em Key Words}\,: 
$q$-hypergeometric series, $(p,q)$-hypergeometric series, $q$-special functions, 
$(p,q)$-special functions, $q$-binomial coefficients, Rogers-Szeg\"{o} polynomial, 
$(p,q)$-binomial coefficients, $(p,q)$-Rogers-Szeg\"{o} polynomial, quantum groups, 
$q$-oscillator, $(p,q)$-oscillator, $(p,q)$-Steiltjes-Wigert polynomial, continuous 
$(p,q)$-Hermite polynomial.}
\footnote{{\em AMS Classification Numbers}\,: 33D45, 33D80, 33D90.} \\
\bigskip
R. Jagannathan\footnote{Formerly of MATSCIENCE, The Institute of Mathematical Sciences, Chennai. \\ 
{\em Email}: {\tt jagan@cmi.ac.in, jagan@imsc.res.in}}\\
\smallskip 
{\em Chennai Mathematical Institute \\ 
Plot H1, SIPCOT IT Park, Padur P.O. \\
Siruseri 603103, Tamilnadu, India} \\
\medskip 
 and \\
\medskip 
R. Sridhar\footnote{Formerly of MATSCIENCE, The Institute of Mathematical Sciences, Chennai. \\ 
{\em Email}: {\tt sridhar@imsc.res.in}}\\
\smallskip 
{\em 30/1, Sundar Enclave, Valmiki Street, Thiruvanmiyur \\  
Chennai 600041, Tamilnadu, India} 
\end{center} 
\medskip 
\centerline{\sf Dedicated to the memory of Professor Alladi Ramakrishnan} 
\vspace{1cm}
\begin{abstract}
A $(p,q)$-analogue of the classical Rogers-Szeg\"{o} polynomial is defined by replacing the 
$q$-binomial coefficient in it by the $(p,q)$-binomial coefficient corresponding to the definition 
of $(p,q)$-number as $[n]_{p,q} = (p^n-q^n)/(p-q)$.  Exactly like the Rogers-Szeg\"{o} polynomial 
is associated with the  $q$-oscillator algebra it is found that the $(p,q)$-Rogers-Szeg\"{o} 
polynomial is associated with the $(p,q)$-oscillator algebra. 
\end{abstract}

\vfill
\noindent
{\sf To appear in}\,: ``{\em The Legacy of Alladi Ramakrishnan in the Mathematical Sciences}''
(K. Alladi, J. Klauder, C. R. Rao, Editors) Springer, 2010.

\newpage 
\section{Introduction} 
\renewcommand\theequation{\thesection.\arabic{equation}}
\setcounter{equation}{0}

\noindent 
The $q$-oscillator algebra plays a central role in the physical applications of quantum groups 
(for a review of quantum groups and their applications see, {\em e.g.}, \cite{CP}-\cite{CD}).  
It was used~\cite{M}-\cite{H} to extend the Jordan-Schwinger realization of the $sl(2)$ algebra 
in terms of harmonic oscillators to the $q$-analogue of the universal enveloping algebra of 
$sl(2)$, namely, $U_q(sl(2))$.  In order to extend this $q$-oscillator realization of 
$U_q(sl(2))$ to the two-parameter quantum algebra $U_{p,q}(gl(2))$  the $(p,q)$-oscillator algebra 
was introduced in \cite{CJ} (see also \cite{JBM,Aet}).  

Heine's $q$-number, or the basic number,  
\beq
[n]_q = \fr{1-q^n}{1-q}, 
\lb{qnumber} 
\eeq
is well known in the mathematics literature.   The $(p,q)$-oscillator necessitated the 
introduction of the $(p,q)$-number, or the twin-basic number, 
\beq
[n]_{p,q} = \fr{p^n-q^n}{p-q}, 
\lb{pqnumber} 
\eeq
a natural generalization of the $q$-number, such that 
\beq
\lim_{p\ra 1} [n]_{p,q} \lra [n]_q. 
\eeq  
With the introduction of this $(p,q)$-number, the essential elements of the $(p,q)$-calculus, 
namely, $(p,q)$-differentiation, $(p,q)$-integration, and the $(p,q)$-exponential, were also 
studied in~\cite{CJ}.  This led to a more detailed study of $(p,q)$-hypergeometric series and 
$(p,q)$-special functions~\cite{BK}-\cite{JS}.  Meanwhile, the $(p,q)$-binomial coefficients, 
$(p,q)$-Stirling numbers, and the $(p,q)$-binomial theorem for noncommutative operators had 
been studied~\cite{KK,SW} in the analysis of certain physical problems.  Interestingly, in the 
same year 1991 the $(p,q)$-number was introduced in the mathematics literature in connection 
with set partition statistics~\cite{WW}.  A very general formalism of deformed hypergeometric 
functions has been developed in~\cite{GGR}.  Some applications of $(p,q)$-hypergeometric series 
in the context of two-parameter quantum groups can be found in~\cite{N,SS}.  

It is noted in~\cite{M} that the classical Rogers-Szeg\"{o} polynomials provide a basis for a 
coordinate representation of the $q$-oscillator.  Several aspects of this close connection 
between the $q$-oscillator algebra and the Rogers-Szeg\"{o} polynomials, and the related 
continuous $q$-Hermite polynomials, have been analysed in detail later~(see, {\em e.g.}, 
\cite{JvJ}-\cite{OS}).   In this paper, after a brief review of the known connection between 
the Rogers-Szeg\"{o} polynomial and the $q$-oscillator, we shall define a $(p,q)$-Rogers-Szeg\"{o} 
polynomial and show that it is connected with the $(p,q)$-oscillator.  

As explained below in section 4, it is not possible to rewrite a $(p,q)$-hypergeometric series, 
or a $(p,q)$-analogue of a $q$-function, as a regular $q$-hypergeometric series or a $q$-function 
routinely by rescaling the independent variable.   Particularly, this is not possible in the case 
of the $(p,q)$-Rogers-Szeg\"{o} polynomial considered here. 

\section{Harmonic oscillator} 
\renewcommand\theequation{\thesection.\arabic{equation}}
\setcounter{equation}{0}

\noindent 
The harmonic oscillator is associated with the creation (or raising) operator $\h{a}_+$, the 
annihilation (or lowering) operator $\h{a}_-$, and the number operator $\h{n}$ satisfying the 
algebra 
\beq 
\lsb \h{n} , \h{a}_+ \rsb = \h{a}_+, \qquad  
\lsb \h{n} , \h{a}_- \rsb = -\h{a}_-, \qquad 
\lsb \h{a}_- , \h{a}_+ \rsb = 1, 
\label{ho}
\eeq
where $\lsb \h{A} , \h{B} \rsb$ stands for the commutator $\h{A}\h{B} - \h{B}\h{A}$.  Note that  
\beq
\h{n} = \h{a}_+\h{a}_-. 
\eeq 
Let 
\beq
h_n(x) = (1+x)^n, \qquad 
\psi_n(x) = \fr{1}{\sqrt{n!}}h_n(x).
\eeq 
It follows that 
\bea 
\fr{d}{dx}\psi_n(x) & = & \sqrt{n}\psi_{n-1}(x), \lb{holo} \\
(1+x)\psi_n(x) & = & \sqrt{n+1}\psi_{n+1}(x),  \lb{horo} \\  
(1+x)\fr{d}{dx}\psi_n(x) & = & n\psi_n(x),  \lb{diffeqn} \\ 
\fr{d}{dx}\lrb (1+x)\psi_n(x) \rrb & = & (n+1)\psi_n(x).  \lb{nplus1} 
\eea
Thus, it is clear that the set $\{\psi_n(x) | n = 0,1,2,\ldots \}$ forms a basis for the following 
Bargmann-Fock realization of the oscillator algebra (\ref{ho}): 
\beq 
\h{a}_+ = (1+x), \qquad 
\h{a}_- = \fr{d}{dx}, \qquad 
\h{n} = (1+x)\fr{d}{dx}. 
\eeq 
It is to be noted that equations (\ref{horo}) and (\ref{diffeqn}) are, respectively, the recurrence 
relation and the differential equation for $\psi_n(x)$.  

Let 
$\lcb \h{a}^{(1)}_-, \h(a)^{(1)}_+, \h{n}^{(1)} \rcb$ and $\lcb \h{a}^{(2)}_-, 
\h{a}^{(2)}_+, \h{n}^{(2)} \rcb$ 
be two sets of oscillator operators  each satisfying the algebra (\ref{ho}) and commuting with 
each other.   Then, for the generators of $sl(2)$ satisfying the Lie algebra, 
\beq 
\lsb x_0 , x_+ \rsb = x_+, \qquad 
\lsb x_0 , x_- \rsb = -x_-, \qquad 
\lsb x_- , x_+ \rsb = 2x_0,  
\eeq 
one has the Jordan-Schwinger realization 
\beq
x_+ = \h{a}^{(1)}_+\h{a}^{(2)}_-, \qquad 
x_ - = \h{a}^{(2)}_+\h{a}^{(1)}_-, \qquad 
x_0 = \fr{1}{2}\lrb \h{n}^{(1)} - \h{n}^{(2)}\rrb. 
\eeq 

\section{$q$-Oscillator and the Rogers-Szeg\"{o} polynomial}  
\renewcommand\theequation{\thesection.\arabic{equation}}
\setcounter{equation}{0}

\noindent 
When $U(sl(2))$, the universal enveloping algebra of $sl(2)$, is $q$-deformed the resulting 
$U_q(sl(2))$ has generators $\lcb X_-, X_+, X_0 \rcb$ satisfying the algebra 
\bea
\lsb X_0 , X_+ \rsb & = & X_+, \qquad 
\lsb X_0 , X_- \rsb = -X_-, \nn \\ 
X_-+X_- - q^{-1}X_-X_+ & = & \fr{1 - q^{2X_0}}{1 - q} = \lsb 2X_0 \rsb_q. 
\lb{uqsl2} 
\eea    
The $q$-oscillator is associated with the annihilation operator $\h{A}_-$, creation operator 
$\h{A}_+$ and the number operator  $\h{N}$ satisfying the algebra 
\beq 
\lsb \h{N} , \h{A}_- \rsb = -\h{A}_-, \qquad  
\lsb \h{N} , \h{A}_+ \rsb = \h{A}_+, \qquad 
\h{A}_-\h{A}_+ - q\h{A}_+\h{A}_- = 1.   
\lb{qo} 
\eeq 
It should be noted that in this case $\h{N} \neq \h{A}_+\h{A}_-$.  Instead we have 
\beq
\h{A}_+\h{A}_- = \fr{1-q^{\h{N}}}{1-q} = \lsb \h{N} \rsb_q,  
\eeq 
and 
\beq
\h{A}_-\h{A}_+ = \fr{1-q^{\h{N}+1}}{1-q} = \lsb \h{N}+1 \rsb_q,  
\eeq    
Now, let 
$\lcb \h{A}^{(1)}_-, \h(A)^{(1)}_+, \h{N}^{(1)} \rcb$ and $\lcb \h{A}^{(2)}_-, 
\h{A}^{(2)}_+, \h{N}^{(2)} \rcb$  
be two sets of $q$-oscillator operators each satisfying the algebra (\ref{qo}) and commuting 
with each other.   Then, taking  
\beq
X_+ = \h{A}^{(1)}_+q^{-\h{N}^{(2)}/2}\h{A}^{(2)}_-, \quad 
X_- = \h{A}^{(2)}_+q^{-\h{N}^{(2)}/2}\h{A}^{(1)}_-,  \quad 
X_0 = \fr{1}{2}\lrb \h{N}^{(1)}-\h{N}^{(2)}\rrb, 
\eeq 
we get a Jordan-Schwinger-type realization of the $U_q(sl(2))$ (\ref{uqsl2}).   

Now, we have to recall some definitions from the theory of $q$-series.  The $q$-shifted factorial 
is defined as 
\beq
(a;q)_n = \lcb \ba{ll} 1, & \mbox{for} \ \ n = 0, \\ 
          \prod_{k=0}^{n-1} (1-aq^k), &  \mbox{for} \ \  n = 1,2,\ldots. \\ \ea \right.   
\eeq 
The $q$-binomial coefficient is defined by 
\beq
\lsb\ba{c} n \\ k \\ \ea \rsb_q = \fr{(q;q)_n}{(q;q)_k (q;q)_{n-k}}.  
\eeq 
For more details on $q$-series see \cite{GR}.  With the definition 
\beq
[0]_q! = 1, \quad 
[n]_q! = [n]_q[n-1]_q\ldots[2]_q[1]_q, \quad 
\mbox{for} \ \  n = 1,2,\ldots, 
\eeq
we have 
\beq
\lsb\ba{c} n \\ k \\ \ea \rsb_q = \fr{[n]_q!}{[k]_q![n-k]_q!}
\eeq 
and 
\beq
\lim_{q \ra 1}  \lsb\ba{c} n \\ k \\ \ea \rsb_q = \fr{n!}{k!(n-k)!} = \lrb\ba{c} n \\ k \\ \ea \rrb. 
\eeq 
The Rogers-Szeg\"{o} polynomial is defined as 
\beq
H_n(x;q) = \sum_{k=0}^n \lsb\ba{c} n \\ k \\ \ea \rsb_q x^k.
\lb{RS} 
\eeq 
This can be naturally expected to be related to the basis of a realization of the $q$-oscillator 
since 
\beq
\lim_{q\ra 1} H_n(x;q) = h_n(x), 
\eeq 
and the $q$-oscillator becomes the ordinary oscillator in the limit $q \lra 1$.  

To exhibit the relation between $H_n(x;q)$ and the $q$-oscillator we shall closely follow \cite{G}, 
though our treatment is slightly different.  Let us define 
\beq
\psi_n(x;q) = \fr{1}{\sqrt{[n]_q!}}H_n(x;q) 
                    =  \fr{1}{\sqrt{[n]_q!}}\sum_{k=0}^n \lsb\ba{c} n \\ k \\ \ea \rsb_q x^k.
\eeq 
The Jackson $q$-difference operator is defined by 
\beq
\h{D}_qf(x) = \fr{f(x)-f(qx)}{(1-q)x}. 
\eeq
It is straightforward to see that 
\beq
\h{D}_q\psi_n(x;q) = \sqrt{[n]_q}\psi_{n-1}(x;q).
\lb{qlo}
\eeq	
The $q$-binomial coefficients obey the recurrence relation 
\beq
\lsb\ba{c} n+1 \\ k \\ \ea \rsb_q = \lsb\ba{c} n \\ 
                  k \\ \ea \rsb_q + \lsb\ba{c} n \\ k-1 \\ \ea \rsb_q 
                  - (1-q^n)\lsb\ba{c} n-1 \\ k-1 \\ \ea \rsb_q. 
\lb{qid} 
\eeq
>From this it follows that $\psi_n(x;q)$ satisfies the recurrence relation 
\beq 
\sqrt{[n+1]_q}\psi_{n+1}(x;q) = (1+x)\psi_n(x;q) - x(1-q)\sqrt{[n]_q}\psi_{n-1}(x;q).
\eeq  
Using (\ref{qlo}) we can write this relation as 
\beq
\lsb (1+x) - (1-q)x\h{D}_q \rsb \psi_n(x;q) = \sqrt{[n+1]_q}\psi_{n+1}(x;q).
\lb{qro} 
\eeq 

Thus, it is seen from (\ref{qlo}) and (\ref{qro}) that the set of polynomials 
$\{ \psi_n(x;q) | n=0,1,2,\ldots\}$ provides a basis for a realization of the $q$-oscillator 
algebra (\ref{qo}) as follows.  Let us define the number operator $\h{N}$ formally as 
\beq 
\h{N}\psi_n(x;q) = n\psi_n(x;q).
\eeq
Note that 
\beq 
\h{D}_q^{n+1}\psi_n(x;q) = 0, 
\eeq 
and $\h{D}_q^m\psi_n(x;q) \neq 0$ for any $m < n+1$.  Then,
\bea  
\h{A}_-\psi_n(x;q) & = & \h{D}_q\psi_n(x;q) = \sqrt{[n]_q}\psi_{n-1}(x;q), \lb{qolo} \\ 
\h{A}_+\psi_n(x;q) & = & \lsb (1+x) - (1-q)x\h{D}_q \rsb \psi_n(x;q) \nn \\  
                                 & = & \sqrt{[n+1]_q}\psi_{n+1}(x;q), \lb{qoro}  \\
\h{A}_+\h{A}_-\psi_n(x;q) & = & [n]_q\psi_n(x;q) = [\h{N}]_q\psi_n(x;q), \lb{qonq} \\ 
\h{A}_-\h{A}_+\psi_n(x;q) & = & [n+1]_q\psi_n(x;q) = [\h{N}+1]_q\psi_n(x;q) \lb{qon1}.   
\eea 
>From this one can easily verify that the relations in (\ref{qo}) are satisfied by 
$\{ \h{A}_+, \h{A}_-, \h{N}\}$.   It may be noted that these relations (\ref{qolo})-(\ref{qon1}) 
are the $q$-generalizations of the harmonic oscillator relations (\ref{holo})-(\ref{nplus1}) to 
which they reduce in the limit $q \lra 1$.   Substituting the explicit expressions for $\h{A}_+$ 
and $\h{A}_-$ in (\ref{qonq}) we get the $q$-differential equation for $\psi_n(x;q)$ 
(or $H_n(x;q)$; see \cite{G})\,:
\beq
\lrb (1-q)x\h{D}_q^2 - (1+x)\h{D}_q + [n]_q \rrb \psi_n(x;q) = 0. 
\lb{rsdiffeqn} 
\eeq 
which reduces to (\ref{diffeqn}) in the limit $q \lra 1$.  

\section{$(p,q)$-Oscillator and the $(p,q)$-Rogers-Szeg\"{o} polynomial}  
\renewcommand\theequation{\thesection.\arabic{equation}}
\setcounter{equation}{0}

\noindent
A genuine two-parameter quantum deformation exists only for $U(gl(2))$ and not for $U(sl(2))$.  
The two-parameter deformation of $U(gl(2))$ leads to $U_{p,q}(gl(2))$ which is generated by 
$\{\hsx_0, \hsx_+, \hsx_-\}$  satisfying the commutation relations 
\bea
\lsb \hsx_0 , \hsx_+ \rsb & = & \hsx_+, \qquad 
\lsb \hsx_0 , \hsx_- \rsb = -\hsx_-, \nn \\
\hsx_+\hsx_- - (pq)^{-1}\hsx_-\hsx_+ & = & \fr{p^{2\hsx_0} - q^{2\hsx_0}}{p-q} = [2\hsx_0]_{p,q}, 
\lb{ugl}
\eea
and a central element $\hsz$ which we shall ignore for the present purpose.  Here, 
$[\phantom{X}]_{p,q}$ is as defined in (\ref{pqnumber}). 

To get an oscillator realization of the algebra (\ref{ugl}) we need the $(p,q)$-oscillator 
algebra defined by 
\beq 
\lsb \hsn , \hsa_+ \rsb = \hsa_+, \quad 
\lsb \hsn , \hsa_- \rsb = -\hsa_-, \quad 
\hsa_-\hsa_+ - q\hsa_+\hsa_- = p^{\hsn}, 
\lb{pqo}
\eeq
where $\{ \hsa_+, \hsa_-, \hsn \}$ are, respectively, the creation, annihilation, and number 
operators.  In this case 
\beq
\hsa_+\hsa_- = \fr{p^{\hsn} - q^{\hsn}}{p-q} = [\hsn]_{p,q}, \qquad 
\hsa_-\hsa_+ = \fr{p^{\hsn+1} - q^{\hsn+1}}{p-q} = [\hsn+1]_{p,q}.  
\eeq
Note the symmetry of this relation under the exchange of $p$ and $q$.  So, the last relation 
in (\ref{pqo}) can also be  taken, equivalently, as  
\beq
\hsa_-\hsa_+ - p\hsa_+\hsa_- = q^{\hsn}. 
\eeq 
The $(p,q)$-oscillator unifies several special cases of $q$-oscillators, including the $q$-fermion 
oscillator \cite{PV}. Now, we shall show that a $(p,q)$-deformation of the Rogers-Szeg\"{o} 
polynomial can be used for a realization of the $(p,q)$-oscillator algebra exactly in the same way 
as the classical Rogers-Szeg\"{o} polynomial is used for a realization of the $q$-oscillator algebra 
as seen above.  To this end we proceed as follows.  

First, let us recall some essential elements of the $(p,q)$-series (for more details see 
\cite{J,JS}).  The $(p,q)$-shifted factorial is defined by 
\bea
(a,b;p,q)_n = \lcb \ba{ll} 1, & \mbox{for} \ \ n = 0, \\ 
              \prod_{k=0}^{n-1} (ap^k-bq^k), &  \mbox{for} \ \  n = 1,2,\ldots. \\ \ea \right.   
\eea 
Note that 
\beq 
(a,b;p,q)_n = a^np^{n(n-1)/2}(b/a;q/p)_n. 
\eeq 
In view of this, it is not possible to rewrite a $(p,q)$-hypergeometric series, or a  
$(p,q)$-analogue of a $q$-function, routinely as a $q/p$-hypergeometric series or a $q/p$-function, 
with the same or a rescaled independent variable, unless the factors depending on $a$ and $p$ in 
the numerator and the denominator cancel in each term, or are such that the uncancelled factor in 
each term is of the same power as the independent variable.  The $(p,q)$-binomial coefficient is 
defined by 
\beq
\lsb\ba{c} n \\ k \\ \ea \rsb_{p,q} = \fr{(p,q;p,q)_n}{(p,q;p,q)_k (p,q;p,q)_{n-k}}.  
\eeq 
With the definition 
\beq
[0]_{p,q}! = 1, \quad 
[n]_{p,q}! = [n]_{p,q}[n-1]_{p,q}\ldots[2]_{p,q}[1]_{p,q}, \quad 
\mbox{for} \ \  n = 1,2,\ldots, 
\eeq
we have 
\beq
\lsb\ba{c} n \\ k \\ \ea \rsb_{p,q} = \fr{[n]_{p,q}!}{[k]_{p,q}![n-k]_{p,q}!} 
\eeq 
and 
\beq
\lim_{p \ra 1}  \lsb\ba{c} n \\ k \\ \ea\rsb_{p,q} = \lsb\ba{c} n \\ k \\ \ea\rsb_q. 
\eeq 

Let us now define the $(p,q)$-Rogers-Szeg\"{o} polynomial as 
\beq
H_n(x;p,q) = \sum_{k=0}^n \lsb\ba{c} n \\ k \\ \ea \rsb_{p,q} x^k, 
\eeq 
and take 
\beq
\psi_n(x;p,q) = \fr{1}{\sqrt{[n]_{p,q}!}}H_n(x;p,q) 
              = \fr{1}{\sqrt{[n]_{p,q}!}}\sum_{k=0}^n \lsb\ba{c} n \\ k \\ \ea \rsb_{p,q} x^k.  
\eeq 
Note that 
\beq
\lsb\ba{c} n \\ k \\ \ea \rsb_{p,q} = p^{k(n-k)}\lsb\ba{c} n \\ k \\ \ea \rsb_{q/p}, 
\eeq 
and the presence of the factor $p^{-k^2}$ makes it impossible to rescale $x$ in any way and 
hence rewrite the $(p,q)$-Rogers-Szeg\"{o} polynomial $H_n(x;p,q)$ as a regular Rogers-Szeg\"{o} 
polynomial (\ref{RS}).  Recalling the definition of the $(p,q)$-difference operator \cite{CJ},   
\beq
\h{D}_{p,q}f(x) = \fr{f(px)-f(qx)}{(p-q)x}, 
\eeq 
it is seen that 
\beq
\h{D}_{p,q}\psi_n(x;p,q) = \sqrt{[n]_{p,q}}\psi_{n-1}(x;p,q). 
\lb{pqlo}
\eeq 
The  $(p,q)$-analogue of (\ref{qid}) is given by 
\beq
\lsb\ba{c} n+1 \\ k \\ \ea \rsb_{p,q} = p^k \lsb\ba{c} n \\ k \\ \ea \rsb_{p,q}
   + p^{n-k+1}\lsb\ba{c} n \\ k-1 \\ \ea \rsb_{p,q} - (p^n-q^n)\lsb\ba{c} n-1 \\ 
                              k-1 \\ \ea \rsb_{p,q}.
\lb{pqid} 
\eeq 
For a detailed study of the $(p,q)$-binomial coefficients see \cite{C}.  This identity (\ref{pqid}) 
leads to the following recurrence relation for $\psi_n(x;p,q)$\,: 
\bea
 \sqrt{[n+1]_{p,q}}\psi_{n+1}(x;p,q) & = & \psi_n(px;p,q) + xp^n\psi_n(p^{-1}x;p,q) \nn \\ 
                    &    & \ \  - x(p-q)\sqrt{[n]_{p,q}}\psi_{n-1}(x;p,q).
\lb{pqro} 
\eea

To obtain a realization of the $(p,q)$-oscillator algebra in the basis provided by 
$\{\psi_n(x;p,q)\,|\,n = 0,1,2,\ldots\}$ let us proceed as follows.  As before, define the 
number operator $\hsn$ formally as 
\beq
\hsn\psi_n(x;p,q) = n\psi_n(x;p,q).  
\lb{pqnumbo} 
\eeq 
Note that 
\beq
\h{D}_{p,q}^{n+1}\psi_n(x;p,q) = 0, 
\eeq 
and $\h{D}_{p,q}^m\psi_n(x;p,q) \neq 0$ for any $m < n+1$.  Then, with the scaling operator 
defined by 
\beq
\he_sf(x) = f(sx),  
\eeq 
it readily follows from (\ref{pqlo}) and (\ref{pqro}) that we can write 
\bea 
\hsa_-\psi_n(x;p,q) & = & 
         \h{D}_{p,q}\psi_n(x;p,q) = \sqrt{[n]_{p,q}}\psi_{n-1}(x;p,q), \lb{pqolo} \\
\hsa_+\psi_n(x;p,q) & = & 
         \lrb \he_p + x\he_{p^{-1}}p^{\hsn}-x(p-q)\h{D}_{p,q} \rrb \psi_n(x;p,q) \nn \\ 
                   & = & \sqrt{[n+1]_{p,q}}\psi_{n+1}(x;p,q), \lb{pqoro} \\ 
\hsa_+\hsa_-\psi_n(x;p,q) & = & 
         [n]_{p,q}\psi_n(x;p,q) = [\hsn]_{p,q}\psi(x;p,q), \lb{pqonpq} \\ 
\hsa_-\hsa_+\psi_n(x;p,q) & = & 
         [n+1]_{p,q}\psi_n(x;p,q) = [\hsn + 1]_{p,q}\psi_n(x;p,q) \lb{pqon1}, 
\eea 
which generalize the corresponding results (\ref{qolo})-(\ref{qon1}) for the $q$-oscillator; 
when $p \lra 1$, (\ref{pqolo})-(\ref{pqon1}) reduce to (\ref{qolo})-(\ref{qon1}).   It is 
straightforward to verify that the realizations of $\{ \hsa_-,\hsa_+,\hsn \}$ in 
(\ref{pqnumbo}), (\ref{pqolo}) and (\ref{pqoro}) satisfy the required relations of the 
$(p,q)$-oscillator algebra (\ref{pqo}).   Using (\ref{pqolo}) and (\ref{pqoro}) in 
(\ref{pqonpq}) we get the $(p,q)$-differential equation satisfied by $\psi_n(x;p,q)$ as 
\beq
\lsb (p-q)x\h{D}_{p,q}^2 - \lrb \he_p+p^{n-1}x\he_{p^{-1}} \rrb
                \h{D}_{p,q} + [n]_{p,q} \rsb \psi_n(x;p,q) = 0,  
\eeq 
which reduces to (\ref{rsdiffeqn}) in the limit $p \lra 1$.   

\section{Conclusion} 
\renewcommand\theequation{\thesection.\arabic{equation}}
\setcounter{equation}{0}

\noindent 
Let us conclude with a few remarks.   \\ 
 
By choosing different values for $p$ and $q$ one can study different special cases of 
$\psi_n(x;p,q)$.  This is particularly important since there exist many versions of 
$q$-oscillators which are special cases of the $(p,q)$-oscillator.   For example, the 
$q$-oscillator originally used in connection with $U_q(su(2))$ \cite{M}-\cite{H}, and more 
popular in physics literature, corresponds to the choice $p = q^{-1}$.  From the above it is 
clear that this oscillator can be realized through the polynomials 
\beq
H_n(x;q^{-1},q) = \sum_{k=0}^n \lsb\ba{c} n \\ k \\ \ea \rsb_{q^2}q^{-k(n-k)}x^k.  
\eeq 
It may be emphasized again that though $H_n(x;q^{-1},q)$ is a function with only a single 
$q$-parameter it cannot be rewritten as a regular Rogers-Szeg\"{o} polynomial.  

In \cite{G} the raising and lowering operators for the Stieltjes-Wigert polynomial have been 
obtained using the fact that this polynomial is just the Rogers-Szeg\"{o} polynomial with $q$ 
replaced by $q^{-1}$ and it has been shown that these raising and lowering operators of the 
Stieltjes-Wigert polynomial provide a realization of the single-parameter deformed oscillator 
with $q$ replaced by $q^{-1}$.   Now, it is clear that one can study the $(p,q)$-Steiltjes-Wigert 
polynomials similarly by replacing $p$ and $q$, respectively, by $p^{-1}$ and $q^{-1}$ in the 
above formalism.  Thus, the $(p,q)$-Steiltjes-Wigert polynomial is given by 
\beq
G_n(x;p,q) = H_n(x;p^{-1},q^{-1}) = 
      \sum_{k=0}^n \lsb\ba{c} n \\ k \\ \ea \rsb_{p,q}(pq)^{-k(n-k)}x^k, 
\eeq 
which becomes the usual Steiltjes-Wigert polynomial in the limit $p \lra 1$.  

The continuous $q$-Hermite polynomial is defined as 
\beq
{\sf H}_n(\cos\theta|q) = e^{-in\theta}H_n(e^{2i\theta};q).
\eeq
It is clear that one can define the continuous $(p,q)$-Hermite polynomial in an analogous way as 
\beq
{\sf H}_n(\cos\theta|p,q) = e^{-in\theta}H_n(e^{2i\theta};p,q).
\eeq
This has already been suggested in \cite{JS} without any further study.  It should be 
worthwhile to study the $(p,q)$-Rogers-Szeg\"{o} polynomial, the $(p,q)$-Steiltjes-Wigert 
polynomial and the continuous $(p,q)$-Hermite polynomial in detail.  

\vspace{1cm} 

\noindent 
{\em Acknowledgement}\,: {\sf We would like to recall, with gratitude and pride, our long 
association with Professor Alladi Ramakrishnan who pioneered research in theoretical physics 
in south India, founded MATSCIENCE, The Institute of Mathematical Sciences, Chennai, and 
directed it for more than two decades.  We are his disciples and our scientific careers were 
moulded in his school.}

\end{document}